\documentclass[12pt]{article}
\usepackage[english]{babel}
\usepackage{amsfonts, amsmath}
\def\c{\frak{c}}
\def\g{\frak{g}}
\def\t{\frak{t}}
\def\p{\frak{p}}
\def\a{\frak{a}}

\def\k{\frak{k}}

\def\n{\frak{n}}

\def\i{{\rm i}}
\def\h{\frak{h}}
\def\c{\frak{c}}
\def\Ad{{\rm Ad}}
\def\ad{{\rm ad}}
\def\beq{\begin{equation}}
\def\C{\mathord{\Bbb C}}
\def\R{\mathord{\Bbb R}}
\def\Z{\mathord{\Bbb Z}}

\def\3{{\ss}}
\def\x{\times}
\def\.{\cdot}

\def\<{\langle}
\def\>{\rangle}
\def\^{\wedge}

\def\ms{\medskip\noindent}
\def\proof{\ms{\bf Proof.}\ }
\def\endproof{\nopagebreak\hfill$\framebox[2mm]{}$\ms}

\author{J.-H. Eschenburg and A.L. Mare}

\title{Steepest descent on real flag manifolds}

\begin{document}


\maketitle

\newtheorem{theorem}{\bf Theorem}[section]
\newtheorem{corollary}[theorem]{\bf Corollary}
\newtheorem{proposition}[theorem]{\bf Proposition}
\newtheorem{lemma}[theorem]{\bf Lemma}
\newtheorem{definition}[theorem]{\bf Definition}

\setcounter{tocdepth}{3}

\section{Introduction}

Among the compact homogeneous spaces, a very distinguished
subclass is formed by the {\it (generalized) real flag
manifolds} which by definition are the orbits of the isotropy
representations of Riemannian symmetric spaces ({\it s-orbits}).
This class contains most compact symmetric spaces (e.g. all
hermitian ones), all classical flag manifolds over real, complex
and quaternionic vector spaces, all adjoint orbits of compact
Lie groups ({\it generalized complex flag manifolds}) and many
others. They form the main examples for isoparametric
submanifolds and their focal manifolds (so called {\it constant
principal curvature manifolds}); in fact for most
codimensions these are the only such spaces (cf.
\cite{hot}, \cite{t}, \cite{o}).

Any real flag manifold $M$ enjoys two very peculiar geometric
properties: It carries a transitive action of a {\it noncompact}
Lie group $G$, and it is embedded in euclidean space as a {\it
taut} submanifold, i.e. almost all height or coordinate functions
are perfect Morse functions (at least for $\Z/2$-coefficients).
The aim of our paper is to link these two properties by the
following theorem: {\sl The gradient flow of any height function is a
one-parameter subgroup of $G$ where the gradient is defined with
respect to a suitable homogeneous metric $s$ on $M$; in the case
where $M$ is an adjoint orbit $s$ is a homogeneous K\"ahler
metric.} In other words, the lines of steepest descend for the
height function (gradient flow lines) are obtained by applying a
one-parameter subgroup of $G$. This is an elementary fact when $M$
is a euclidean sphere  and $G$ its conformal group: The gradient
of any height function is a conformal vector field. For adjoint
orbits, the {\it (generalized) complex flag manifolds}, this fact
was observed earlier by Guest and Ohnita \cite{go}. Our more general
result can be derived from their theorem since real flag manifolds
are contained in complex flag manifolds as fixed point sets  under
certain involutions. However a short direct proof might be
desirable. We would like to thank Martin Guest and Peter Quast 
for hints and discussion.

\section{Root space decomposition}

Let $P = G/K$ be a symmetric space of noncompact type where $G$ is a connected
noncompact semisimple Lie group and $K\subset G$ a maximal compact subgroup
(cf. \cite{h}). Let $\sigma$ be the corresponding involution on $G$ with
fixed point set $K$. Consider the corresponding  Cartan decomposition
\beq
    \g=\k \oplus \p
\end{equation}
where $\g$ and $\k$ denote the Lie algebras of $G$ and $K$ and
$\p$ the $(-1)$-eigenspace of (the differential of) $\sigma$. The
adjoint action of $K$ leaves $\p$ invariant; this is the isotropy
representation of $P$. As usual we consider an $\Ad(G)$-invariant
indefinite inner product $b$ on $\g$ with $b > 0$ on $\p$ and $b <
0$ on $\k$ (e.g. the Killing form) and define a positive definite
inner product $\<\ ,\ \>$ on $\g$ which is $b$ on $\p$ and $-b$ on
$\k$, and such that $\p \perp \k$. Then any $\ad(x)$ with $x \in
\p$ is self adjoint on $\g$. Hence a maximal abelian subspace $\a
\subset \p$ gives rise to a family of mutually commuting
selfadjoint endomorphisms $\ad(x)$ of $\g$ where $x\in \a$. These
have a common eigenspace decomposition \beq
    \g=\sum_{\alpha\in \hat R}\g_{\alpha},\ \
    \g_{\alpha}:=\{z\in \g;\ [x,z]=\alpha(x)z\ \forall x\in \a\}
\end{equation}
where $\hat R\subset \a^*$ is the set of roots including $0$.

After fixing some arbitrary $x \in \a$, there are three disjoint subsets of
$\hat R$, formed by the roots $\alpha$ with $\alpha(x) > 0$, $\alpha(x) = 0$
and $\alpha(x) < 0$, respectively. Hence we have the decomposition
\beq \label{decompo}
    \g = \n_+ \oplus \c \oplus \n_-
\end{equation}
with
\beq
    \n_+ = \sum_{\alpha(x)>0} \g_\alpha,\ \
    \c   = \sum_{\alpha(x)=0} \g_\alpha,\ \
    \n_- = \sum_{\alpha(x)<0} \g_\alpha.
\end{equation}
Since $\sigma = -id$ on $\a$, we have
\beq
    (\g_\alpha)^\sigma = \g_{-\alpha},
\end{equation}
hence $\sigma$ interchanges $\n_+$ and $\n_-$. Any $a \in \g$ allows a unique
decomposition $a = a_- + a_o + a_+$ with $a_\pm \in \n_\pm$ and $a_o \in \c$,
and if $a \in \k$ or $a \in \p$, we have $a_- = a_+^\sigma$ or $a_- =
-a_+^\sigma$, respectively.

\section{Generalized real flag manifolds}

By definition, generalized real flag manifolds are the orbits of
the isotropy representation of a symmetric space $P = G/K$ of
noncompact type. We consider an isotropy orbit $M \subset \p$. 
For any fixed $x \in M$ we have $M = \Ad(K)x$. The stabilizer subgroup
is 
\beq
    S=\{k\in K;\ \Ad(k)x=x\},
\end{equation}
thus $M$ may be identified with the coset space $K/S$ by $kS \mapsto
{\rm Ad}(k)x$.

Now choose a maximal abelian subalgebra $\a \subset \p$
containing $x$. Such $\a$ is uniquely determined up to applying
${\rm Ad}(s)$ with $s \in S$; this fact is just the conjugacy of 
maximal flat subspaces in the symmetric subspace $\tilde P = 
C/S \subset P$ where $C$ (with Lie algebra $\c$)
is the centralizer of $x$ in $G$.\footnote
	{Geometrically, $\tilde P$ is the union of all geodesics
	parallel to the geodesic $\exp(\R x) \subset P$ which is
	obviously invariant under geodesic reflection, hence it is
	a symmetric subspace.}
Note that each $\Ad(s)$, $s \in S$, commutes with $ad(x)$ and thus 
preserves the eigenspaces of $ad(x)$, hence $\n_- = 
\sum_{\alpha(x)<0} \g_\alpha$ is also invariant under $\Ad(s)$.

Next we show that the natural $K$-action on $M$ can be extended to
a $G$-action. Consider 
\beq
    H = \{g\in G;\ {\rm Ad}(g)(x +\n_-) =x +\n_-\}
\end{equation}
which is a closed subgroup of $G$ containing $S$ as a subgroup.
One can see that the Lie algebra of $H$ is $\h = \c + \n_-$. Note
that $H$ does not depend on the choice of $\a$.
Indeed, if instead of $\a$ we start our construction with
$\a'={\rm Ad}(s)\a$, where $s\in S$, then $\n'_-={\rm Ad}(s)\n_-
= \n_-$,
and we end up with $H'=H$.

\begin{lemma}\label{lema1} 
$K$ acts transitively on the coset space $G/H$ with stabilizer
$S$. Hence $G/H$ can be identified with $K/S = M$.
\end{lemma}

\proof Let $M' \subset G/H$ be the orbit of $eH \in G/H$ under the
subgroup $K \subset G$. We show first that it is open in $G/H$. To
see this it suffices to show  $\k+\h=\g$, i.e. one needs to show
$\g_\alpha \subset \k+\h$ for each $\alpha \in R$ with $\alpha(x)
> 0$. In fact, take $z\in \g_{\alpha}$, decompose it as $z=v+u$
with $v\in \p$ and $u\in \k$ and notice that $z^\sigma = -v+u \in
\g_{-\alpha}\subset \h$. Thus the $K$-orbit $M'$ is open. But it
is also closed in $G/H$ since $K$ is compact. So $M'$ coincides
with $G/H$.

The $K$-stabilizer of $eH \in G/H$ is $K\cap H$; we have to show
$K\cap H = S$. Clearly $S \subset K\cap H$. Vice versa, if $k\in
K\cap H$, then 
	$${\rm Ad}(k)x -x\in \p \cap \n_-.$$ 
But $\p\cap\n_-=0$, because if $z$  belongs to this intersection, then
$z^{\sigma}=-z\in \n_+$, hence $z\in \n_-\cap \n_+=0$. We deduce
that ${\rm Ad}(k)x=x$, which means $k\in S$.

\endproof

\ms

Thus the action of $G$ on $G/H = K/S = M$ is an extension of the
$K$-action.\footnote
    {There is a more geometric description of this action: Consider $\p$ as
    the tangent space of $P = G/K$ at some base point $o \in P$. We may
    project any nonzero $x \in T_oP$ to the infinite boundary $P(\infty)$
    (mind that $P$ is a simply connected space of nonpositive curvature) by
    the map $\pi_\infty(x) = \gamma_x(-\infty)$ where $\gamma_x$ is the geodesic in
    $P$ starting at $o$ with initial vector $x$. The isometry group $G$ acts
    on $P(\infty)$ and leaves $\pi_\infty(M)$ invariant; this is the
    $G$-action. See \cite{bgs} for details.}
We will denote this action $G\x M \to M$ by $(g,x) \mapsto g.x$. When restricted
to $k \in K \subset G$ we have $k.x = \Ad(k)x$. Similarly, the infinitesimal
action $\g \x M \to TM$ will be denoted by $(a,x) \mapsto a.x := {d\over
dt}|_{_{t=0}}\exp(ta).x$ where $a \in \g$, and $a.x =[a,x]$ whenever $a \in \k$.

\section{Steepest descend}

\begin{theorem}\label{main} Let $G/K$ be a symmetric space of noncompact type, 
$\g = \k +\p$ the corresponding Cartan decomposition and $M \subset \p$ a real
flag manifold (isotropy orbit). Let $q \in \p$ and $f : M \to \R$, $f(x) =
\<q,x\>$. Then there is a $K$-invariant Riemannian metric $s : TM \to T^*M$ on
$M$ such that the flow lines $x(t)$ of the $s$-gradient $\nabla^s f = s^{-1}df$
satisfy 
\beq  x(t) = \exp(-tq).x(0) \end{equation}  
\end{theorem}

\proof Choose an arbitrary $x \in M$ and consider the decomposition
(\ref{decompo}) corresponding to $x$. We have to show that $\nabla^sf(x) =
-q.x$. To compute $q.x$ we look for $r \in \k$ with $q.x = r.x$, i.e. $q-r \in
\h$. We have $q = q_o + \sum_{\alpha(x)>0} q_\alpha$ with $q_o \in \c$ and
$q_\alpha = z_\alpha - z_\alpha^\sigma$ for some $z_\alpha \in \g_\alpha$. We
may assume $q_o = 0$ since $q_o.x = 0$. Then we put $r = \sum_+ (z_\alpha +
z_\alpha^\sigma) \in \k$ and hence $r-q = 2\sum_+ z_\alpha^\sigma \in \h$ where
$\sum_+$ always denotes $\sum_{\alpha(x)>0}$. Now 
$$ 
	r.x = [r,x] = -ad(x)\sum_+(z_\alpha + z_\alpha^\sigma) 
		= -\sum_+ \alpha(x)(z_\alpha - z_\alpha^\sigma) 
		= -\sum_+ \alpha(x) q_\alpha. 
$$ 
Any $v \in T_xM$ has a decomposition $v = \sum_+ v_\alpha$ with $v_\alpha \in
(\g_\alpha + \g_{-\alpha})\cap \p$ for $\alpha(x)>0$, and our metric $s$ on
$T_xM$ will be of the form $$ \<v,w\>_s = \sum_+ s_\alpha \<v_\alpha,w_\alpha\>
$$ for certain numbers $s_\alpha > 0$. We have to choose $s$ such that for all
$v \in T_xM$ $$ \<\nabla^sf(x),v\>_s = -\<r.x,v\>_s . $$ The left hand side is
$$ \<\nabla^sf(x),v\>_s = df_xv = \<q,v\> = \sum_+ \<q_\alpha,v_\alpha\>, $$
while the right hand side is $$ -\<r.x,v\>_s = \sum_+ s_\alpha \alpha(x)
\<q_\alpha,v_\alpha\>. $$ Hence we obtain the result by putting 
\beq \label{s}
	s_\alpha = 1/\alpha(x).
\end{equation} \endproof

\section{Extrinsic symmetric spaces}

An {\it extrinsic symmetric space} is a submanifold $M$ in euclidian
space such that $M$ is preserved by the reflections at all of its
(affine) normal spaces. By a result of Ferus \cite{fe} (also cf.
\cite{eh}), after splitting off euclidean factors, $M$ has precisely 
the form of an s-orbit $M = \Ad(K)x_o \subset \p$ where $\p$ corresponds
to a symmetric space $P = G/K$ and where $x_o \in \p$ satisfies
$$\alpha(x_o)\in\{-1, 0, 1\},$$ for all $\alpha\in R$.
In this case the metric $s$ of Theorem \ref{main} agrees to the given inner
product $\<\ ,\ \>$ (cf. (\ref{s})). By applying
Theorem \ref{main} one obtains:

\begin{theorem}
{\it If $M={\rm Ad}(K)x_o\subset \p$ is extrinsic symmetric, then
the gradient lines of the height function $h(x)=\langle
q,x\rangle$ with respect to the metric $\langle \ , \ \rangle$ are
of the form
$$x(t)=\exp(-tq).x(0).$$}
\end{theorem}

\section{Adjoint orbits}

In the particular case of {\it complex} flag manifolds (i.e.
adjoint orbits of compact Lie groups) we will establish relations
between Theorem \ref{main}  and previously known results.

Let $K$ be a compact semisimple Lie group of Lie algebra $\k$, and
$T\subset K$ a maximal torus of Lie algebra $\t$.  Consider the
adjoint orbit $M={\rm Ad}(K)x$ for $x\in \k$.  If $G=K^{\C}$
is the complexification of $K$, then $G/K$ is a non-compact
symmetric space and
	$$\g=\k+\i\k$$
is a Cartan decomposition of $\g={\rm Lie}(G)=\k\otimes \C$ (the
involution $\sigma$ is just the complex conjugation). Since $M$ is
(up to a multiple of $\i$) an isotropy orbit of $G/K$, the results of the
previous section can be applied here, too. The goal of this
section is to point out that the metric
on $M$ for which the  lines of steepest descent are orbits of
one-parameter subgroups of $G$ is well known: It is the K\"ahler 
metric (cf. \cite{go}, \cite{fe}).

It is well known that any adjoint orbit $M = \Ad(K)x$ is a complex
manifold. In fact, in the language of Section 4 we have $M = G/H$
but here $G$ and $H$ are {\it complex} Lie groups and hence $M$ 
is a complex manifold. The corresponding complex structure $J$
on $T_xM$ can be described as follows. Choose a maximal abelian
subalgebra $\t \subset \k$ with $x \in \t$. The corresponding roots 
are considered as real linear forms $\alpha \in \t^*$ while the
eigenvalues of $ad(x)$ are purely imaginary, $\i\alpha(x)$. Let
$\k_\alpha \subset \k \otimes \C$ be the root spaces. Then
$$
	T_xM = \sum_{\alpha(x)>0} \k_\alpha^r
$$
where $\k_\alpha^r = (\k_\alpha + \k_{-\alpha})\cap \k$ is the
{\it real root space}. Now $J$ leaves invariant each $\k_\alpha^r$ and
on $\k_\alpha^r$ it is a multiple of $\ad(x)$:
\beq \label{ad}
	\ad(x) = \alpha(x)J
\end{equation}

The second ingredient for the K\"ahler metric is the K\"ahler form
$\omega$ which is defined as follows: If $v = \ad(a)x$ and $w = \ad(b)x$
are tangent vectors of $M$ at the point $x$, then 
\beq \label{omega}
	\omega_x(v,w) := \<x,[a,b]\> = \<[x,a],b\>
                 = \<v,b\> = \<v,\ad(x)^{-1}w\>
\end{equation}
where $\<\ ,\ \>$ denotes an $Ad(K)$-invariant inner product on $\k$.

Now the K\"ahler metric $(\ ,\ )$ on $T_xM$ is defined as follows. For
$v,w \in T_xM$ we have
\beq \label{kaehler}
	(v,w) = \omega_x(v,Jw).
\end{equation}
Hence from (\ref{ad}) and (\ref{omega}) we obtain
\beq 
	(v,w) = \omega_x(v,Jw) = \<v,\ad(x)^{-1}Jw\> 
	= {1\over\alpha(x)}\<v,w\> = \<v,w\>_s
\end{equation}
(cf. (\ref{s})).
Since $\p = \i\k$ in the present case, we obtain from Theorem \ref{main}: 

\begin{theorem} (cf. \cite{go}, \cite{fe}) Let $K$ be a compact Lie
group and $M = \Ad(K)x_o$ $\subset \k$ an adjoint orbit, equipped with
its K\"ahler metric $(\ ,\ )$ as in (\ref{kaehler}) and acted on by the
complexified group $G = K^{\C}$ as described above.  Let $\<\ ,\ \>$
be the corresponding $\Ad(K)$-invariant inner product on $\k$. Then for
any $q \in \k$ the gradient lines $x(t)$ of the function $f : M \to
\R$, $f(x) = \langle q,x\rangle$, are orbits of a 1-parameter subgroup
of $G = K^{\C}$, namely $$x(t) = \exp(\i t q).x(0).$$  \end{theorem}

\end{document}